\documentclass[12pt]{article}
\usepackage{amssymb, amsmath}

\newtheorem{theorem}{Theorem}[section]

\newcommand{\vare}{\varepsilon}
\newcommand{\n}{\nonumber}

\renewcommand{\a}{\alpha}
\renewcommand{\o}{\omega}

\newcommand{\bb}{\begin{equation}}
\newcommand{\ee}{\end{equation}}
\newcommand{\bq}{\begin{eqnarray}}
\newcommand{\eq}{\end{eqnarray}}
\newcommand{\bqn}{\begin{eqnarray*}}
\newcommand{\eqn}{\end{eqnarray*}}

\begin{document}
\title{Global regularity for a model Navier-Stokes equations on $\Bbb R^3$}
\author{ Dongho Chae\\
\ \\
Department of Mathematics\\
 Chung-Ang
University\\
 Seoul 156-756, KOREA\\
 e-mail: {\it dchae@cau.ac.kr }}
 \date{}
\maketitle
\begin{abstract}
We study a nonlinear parabolic system for a time dependent solenoidal vector field on $\Bbb R^3$.
The nonlinear term of this new model equations  is obtained slightly modifying that of the Navier-Stokes equations. 
The system has the same scaling property and the Galileian  invariance as the Navier-Stokes equations. For such system we prove the global regularity
for a smooth initial data.\end{abstract}
\section{Introduction}
 \setcounter{equation}{0}
 We consider the incompressible Navier-Stokes equations in $\Bbb R^3$,
$$
  (NS)\left\{\aligned  & v_t +v\cdot \nabla v=-\nabla p+\Delta v, \quad (x,t)\in \Bbb R^3\times [0, \infty)\\
     & \nabla \cdot v=0, \quad (x,t)\in \Bbb R^3\times [0, \infty)\\
& v(x,0)=v_0(x),  \quad x\in \Bbb R^3\\
    \endaligned \right.
$$
 where $v=(v_1, v_2, v_3), v_j =v_j (x,t), j=1,2,3,$ is the velocity field, and 
 $p=p(x,t)$ is the pressure. For simplicity we consider the case of zero external force.
 The problem of global regularity/finite time singularity  of solutions to (NS) for a smooth initial data  $v_0$ is an  outstanding open problem.  We know the local in time well-posedness 
 for smooth initial data(\cite{lic, kat}), the global existence of weak solutions(\cite{ler, hop}),  
 and the partial regularity of suitable weak solutions(\cite{sch, caf}). For comprehensive studies 
 of the Cauchy problem of (NS) we refer \cite{lad, lem}.
 Using the vector identity
 $$ v\cdot \nabla v=-v\times \mathrm{curl }\, v+\frac12 \nabla |v|^2,
 $$
One can rewrite the system (NS) in  an equivalent form,
 $$
  (NS)_1\left\{\aligned  & v_t -v\times \o =-\nabla (p+\frac12 |v|^2)+\Delta v,\\
   &\nabla \cdot v=0,\,\o=\nabla \times\ v\\
    & v(x,0)=v_0(x).
    \endaligned \right.
$$
 For further discussion let us recall the definition of the Riesz transform(\cite{ste}) in $\Bbb R^n$, $n\in \Bbb N$. For $f\in L^p (\Bbb R^n)$, $p\in [1, \infty)$, the Riesz transform $R_j$ 
 of $f$ is given by 
 $$
 R_j (f)(x)= c_n \lim_{\vare\to 0}\int_{|y|>\vare} \frac{y_j}{|y|^{n+1}} f(x-y) dy, \quad j=1, \cdots , n,
 $$
 where 
 $$ c_n=\frac{\Gamma \left(\frac{n+1}{2}\right)}{\pi^{(n+1)/2}}.$$
 Let $\hat{f}(\xi)= \int_{\Bbb R^n} e^{2\pi i x\cdot \xi } f(x) dx$ be the Fourier transform of $f$.
Then, the Riesz transform is more conveniently defined by the Fourier transform as
\bb\label{rie}
(\widehat{R_j f})(\xi )=i \frac{\xi_j}{|\xi|} \hat{f}(\xi), \quad i=\sqrt{-1}.
\ee
Let $R_1, R_2, R_3$ be the Riesz transforms in $\Bbb R^3$, and $u=(u_1,u_2, u_3)$ be a vector field on $\Bbb R^3$. Then, we define for scalar function $f$,
$ R(f):=(R_1(f), R_2 (f), R_3 (f))$, and for a vector field $u=(u_1, u_2, u_3)$,
$$R\cdot u:= R_1 (u_1)+R_2 (u_2)+R_3(u_3),
$$
and 
$$ R\times u:= (R_2u_3-R_3u_2, R_3u_1-R_1u_3, R_1u_2-R_2u_1).
$$ 
Then,  using the fact  $ R_1^2 +R_2 ^2 +R_3 ^2 =-I$,  which follows immediately from (\ref{rie}), we have for a vector field $F$, satisfying div $F=0$,
 $$R\times R\times F= R^2F-R(R\cdot F)=-F .$$
Since 
 $$ R\times R\times \{\nabla (p+\frac12 |v|^2)\} =0, \quad  R\times R\times \{ v_t -\Delta v\}=-v_t+\Delta v
 $$
 for $v$ with div $v=0$,  the system $(NS)_1$ can be written as
the following equivalent system.
$$
  (NS)_2\left\{\aligned  & v_t +R\times R\times(v\times \o )=\Delta v,\\
    & v(x,0)=v_0(x), \quad \nabla \cdot v_0=0.
    \endaligned \right.
$$
The implication that $(NS)_2$ follows from $(NS)_1$ is obvious by taking $ R\times R\times$ on $(NS)_1$. For the reverse direction we note that the first equation of $(NS)_2$ can be written as
$$  R\times R\times( v_t -v\times \o -\Delta v)=0, $$
from which it follows that there exists a scalar function $p=p(x,t)$ such that 
$$v_t -v\times \o -\Delta v=-\nabla (p+\frac12 |v|^2).$$
 Note that the divergence free condition of the initial data is preserved by the equations 
 in $(NS)_2$ for smooth solutions, and the solution satisfies div $v=0$ automatically.
 We also observe that in the formulation $(NS)_2$ there exists no pressure term,
 although the nonlocality is now moved to the nonlinear term via the Riesz transforms.
 
We expect that the problem of the global regularity/finite time singularity for the system $(NS)_2$ has the similar difficulty as the original Navier-Stokes equations $(NS)$.
Instead of $(NS)_2$  we consider its modified version:
$$
  (mNS)\left\{\aligned  & v_t +R\times (v\times \o )=\Delta v,\\
    & v(x,0)=v_0(x), \quad \nabla \cdot v_0=0,
    \endaligned \right.
$$
which is obtained by omitting one nonlocal operation $R\times $ in the nonlinear term in the system $(NS)_2$. 
Another interpretation of $(mNS)$ is its relation with the  following Hall equations, which is obtained from the Hall magnetohydrodynamics(Hall-MHD) equations by setting the velocity$=0$, and which represent the major difficult part of the whole system(see \cite{cha1,cha2, cha3, cha4} for more detailed studies of the Hall-MHD system).
$$
  (\mathrm{Hall})\left\{\aligned  & B_t +\nabla \times (B\times (\nabla \times B ))=\Delta B,\quad (x,t)\in \Bbb R^3\times (0, \infty)\\
    & B(x,0)=B_0(x), \quad \nabla \cdot B_0=0,\quad x\in \Bbb R^3
    \endaligned \right.
$$
where $B=(B_1, B_2, B_3)$, $B_j=B_j (x,t), j=1,2,3$ is the magnetic field.
Comparing with $(\mathrm{Hall})$, we find that the nonlinear term of $(mNS)$ is regularized by one derivative,  in the sense that $\nabla \times$ in front of the nonlinear term (Hall) is replaced by  $ R\times =\nabla (-\Delta )^{-1/2}\times$
in (mMS).
We note that the symmetry properties(say, the Galilean invariance and the scaling symmetry) of  $(mNS)$
are the same as the original Navier-Stokes equations. Our  result in the following theorem shows  that for such modified system $(mNS)$ we could show the global regularity 
for a given smooth  initial data.  
\begin{theorem} Let $v_0 \in H^m (\Bbb R^3)$ with $m>5/2$. Then, for all $T\in (0, \infty)$ there exists a unique solution $v \in C([0, T); H^m (\Bbb R^3)) \cap L^2 (0, T; H^{m+1} (\Bbb R^3))$ to the system $(mNS)$.
Moreover, the solution satisfies the inequality:
\bq
\lefteqn{\sup_{0<t<T}\| v(t)\|_{H^m} ^2 +\int_0 ^T \|Dv(s)\|_{H^m} ^2 ds}\n \\
&&\leq \| v_0\|_{H^m} ^2\exp \left\{T \|v_0\|_{L^2}^2 \exp (C\|\Lambda ^{\frac12 }
 v_0\|_{L^2}^2)\right\}\times\n \\
  &&\times\exp\left[ \|\Lambda v_0\|_{L^2} ^2\exp \left\{CT\|v_0\|_{L^2}^2 \exp (C\|\Lambda ^{\frac12 }
 v_0\|_{L^2}^2)  \|\Lambda v_0\|_{L^2} ^2\exp (C\|\Lambda ^{\frac12 }
 v_0\|_{L^2}^2)\right\}\right]\n \\
\eq
for all $T>0$.
\end{theorem}
{\em Remark 1.1 } The result shows that the simple estimates for the nonlinear term 
of $(NS)_2$ is not enough to deduce correct result for the problem of the global regularity/finite  time singularity. We mention here that there are also studies of the other model equations  of the Navier-Stokes system, where the authors show the finite time blow-up (see e.g. \cite{che, fri, gal, mon, ple, tao}).
\section{Proof of Theorem 1.1}
 \setcounter{equation}{0}
 Let us set  the
multi-index $\alpha :=(\alpha_1 , \alpha_2 , \cdots , \alpha_n )\in
(\Bbb
N \cup \{ 0\} )^n $ with $|\alpha |=\alpha_1 +\alpha_2 +\cdots
 +\alpha_n$. Then, $D^\alpha :=D^{\alpha_1}_1 D^{\alpha_2}_2 \cdots
 D^{\alpha_n}_n$, where $D_j =\partial/\partial x_j$, $j=1,2,\cdots,
 n$.
Given $m\in \Bbb N\cup \{0\}$  the Sobolev space,
 $H^m (\Bbb R^n )$ is the Hilbert space of functions
  consisting of functions $f\in L^2 (\Bbb R^n )$
 such that
 $$\|f\|_{H^m}:=\left(\sum_{|\a|\leq m}\int_{\Bbb R^n}|D^\alpha f(x)|^2 dx
 \right)^{\frac{1}{2}} <\infty,$$
  where the derivatives are in the sense of distributions. Given $s\in \Bbb R$, we use the notation $\Lambda^s (f)=(-\Delta )^{\frac{s}{2}} f$, which is
defined by its Fourier transform as
$$
\widehat{\Lambda^s (f)} (\xi)= |\xi|^s \hat{f}(\xi).
$$
 Then, we observe that 
 $$ \widehat{R_j (f)}= i \frac{\xi_j}{|\xi|} \hat{f}=\widehat{\partial_j \Lambda ^{-1}f} (\xi).
 $$
 We use the energy method for the proof of Theorem 1.1(see e.g. \cite{maj}).\\
 \ \\
\noindent{\bf Proof  of Theorem 1.1 } Since the local well-posedness in the Sobolev space $H^m(\Bbb R^3)$ for $m>5/2$ is 
standard, we proceed directly to the global in time a priori estimate.
We take $L^2$ inner product $(mNS)$ by $\Lambda v$ to
obtain
 \bqn
 \lefteqn{\frac12 \frac{d}{dt} \|\Lambda ^{\frac12} v\|_{L^2} ^2
 +\|\Lambda ^{\frac32} v\|_{L^2}^2 =-\int_{\Bbb R^3} R\times (v\times
 \o) \cdot \Lambda v \, dx} \\
 && =-\int_{\Bbb R^3}  (v\times
 \o) \cdot \Lambda R\times v \, dx=-\int_{\Bbb R^3}  (v\times
 \o) \cdot \o \, dx=0.
 \eqn
We have therefore
  \bb\label{en}
 \frac12\|\Lambda ^{\frac12} v(t)\|_{L^2} ^2+\int_0 ^t\|\Lambda ^{\frac32}
 v(s)\|_{L^2}^2ds = \frac12\|\Lambda ^{\frac12} v_0\|_{L^2} ^2
\ee
 Next we take $L^2$ inner product $(mNS)$ by $ v$ to
deduce
 \bqn
 \lefteqn{\frac12 \frac{d}{dt} \| v\|_{L^2} ^2
 +\|\nabla v\|_{L^2}^2 =-\int_{\Bbb R^3} v\times
 \o \cdot R\times  v \, dx} \\
 && \leq \int_{\Bbb R^3} |v||\o||Rv|\, dx\leq
 \|v\|_{L^6}\|\o\|_{L^3} \|Rv\|_{L^2}\\
 &&\leq C\|\nabla v\|_{L^2} \|\Lambda ^{\frac32} v\|\|v\|_{L^2}\\
 &&\leq\frac12 \|\nabla v\|_{L^2}^2 +C\|\Lambda ^{\frac32}
 v\|^2\|v\|_{L^2}^2.
 \eqn
 Hence, we have
\bq\label{est1}
 \| v(t)\|_{L^2} ^2+\int_0 ^t \|\nabla v\|_{L^2}^2ds &\leq&
 \|v_0\|_{L^2}^2 \exp \left(C\int_0 ^t\|\Lambda ^{\frac32}
 v\|^2ds\right)\n \\
 &\leq& \|v_0\|_{L^2}^2 \exp (C\|\Lambda ^{\frac12 }
 v_0\|_{L^2}^2),
 \eq
 where we used the estimate (\ref{en}).
  We now take $L^2$ inner product $(mNS)$ by $\Delta v$ to
have
 \bqn
 \lefteqn{\frac12 \frac{d}{dt} \|\Lambda v\|_{L^2} ^2
 +\|\Delta v\|_{L^2}^2 =-\int_{\Bbb R^3} v\times
 \o \cdot R\times  \Delta v \, dx} \\
&&\leq \|v\|_{L^6} \|\o\|_{L^3}
 \|R\Delta v\|_{L^2}\leq C \|\nabla v\|_{L^2} \|\Lambda ^{\frac32} v\|_{L^2} \|\Delta
 v\|_{L^2}\n \\
 &&\leq
 \frac12\|\Delta
 v\|_{L^2}^2 + C \|\Lambda ^{\frac32} v\|_{L^2}^2\|\nabla
 v\|_{L^2}^2,
 \eqn
from which we obtain \bq\label{est2}
 \|\Lambda v(t)\|_{L^2} ^2 +\int_0 ^t\|\Delta v\|_{L^2}^2 ds
 &\leq&\|\Lambda v_0\|_{L^2} ^2\exp \left(C\int_0 ^t \|\Lambda ^{\frac32}
 v\|_{L^2}^2ds\right)\n \\
 &\leq&\|\Lambda v_0\|_{L^2} ^2\exp (C\|\Lambda ^{\frac12 }
 v_0\|_{L^2}^2)
 \eq
 by (\ref{en}).
We  operate $D^2$ on $(mNS)$, and take $L^2$ inner product of it by
$D^2 v$ to obtain
 \bqn
 \lefteqn{\frac12 \frac{d}{dt} \|D^2 v\|_{L^2} ^2
 +\|D^3 v\|_{L^2}^2 =-\int_{\Bbb R^3}D( R\times (v\times
 \o)) \cdot D^3 v \, dx} \\
 && =-\int_{\Bbb R^3} (Dv\times
 \o+ v\times D\o ) \cdot R\times D^3  v \, dx\n \\
 &&\leq (\|Dv\|_{L^3}\|\o\|_{L^6} + \|v\|_{L^\infty} \|D\o\|_{L^2}) \|D^3 v\|_{L^2}\n \\
 && \leq C( \|\Lambda ^{\frac32}v\|_{L^2} +\|v\|_{L^\infty})\|D^2
 v\|_{L^2} \|D^3 v\|_{L^2}\n \\
 &&\leq \frac 12\|D^3 v\|_{L^2}^2  +C( \|\Lambda ^{\frac32}v\|_{L^2}^2 +\|v\|_{L^\infty}^2)\|D^2
 v\|_{L^2}^2.
 \eqn
 Hence,
\bq\label{est3}
 \lefteqn{\|D^2 v(t)\|_{L^2} ^2 +\int_0 ^t\|D^3 v\|_{L^2}^2 ds
 \leq\|D^2 v_0\|_{L^2} ^2\exp \left\{C\int_0 ^t (\|\Lambda ^{\frac32}
 v\|_{L^2}^2 +\|v\|_{L^\infty} ^2 )ds\right\}}\n \\
 &&\leq\|\Lambda v_0\|_{L^2} ^2\exp \left\{C \int_0 ^t ( \|v(s)\|_{L^2} ^2 + \|\Delta v(s)\|_{L^2}^2 )ds\right\}.\n \\
 &&\leq \|\Lambda v_0\|_{L^2} ^2\exp \left\{Ct\|v_0\|_{L^2}^2 \exp (C\|\Lambda ^{\frac12 }
 v_0\|_{L^2}^2)  \|\Lambda v_0\|_{L^2} ^2\exp (C\|\Lambda ^{\frac12 }
 v_0\|_{L^2}^2)\right\},\n \\
 \eq
 where we used the estimates (\ref{est1}) and (\ref{est2}).
 Let $m> 5/2$.
 Operating $D^\alpha$ on $(mNS)$ and taking $L^2$ inner product of it by $D^\alpha v$, and summing over
 $|\alpha|\leq m$,  one has
 \bqn
 \lefteqn{\frac12 \frac{d}{dt} \| v\|_{H^m} ^2
 +\|D v\|_{H^m}^2 =-\sum_{|\alpha|\leq m} \int_{\Bbb R^3}D^\alpha  (v\times
 \o) \cdot R\times D^\alpha  v \, dx} \\
 &&=\sum_{|\beta|= |\a|-1, 1\leq |\a|\leq m} \int_{\Bbb R^3}D^\beta  (v\times
 \o) \cdot R\times D^\alpha D  v \, dx - \int_{\Bbb R^3}  (v\times
 \o) \cdot R\times   v \, dx \n \\
  && \leq C\|v\times \o\|_{H^{m-1}} \|R Dv\|_{H^m} \leq  C(\|v\|_{L^\infty} +\|\nabla v\|_{L^\infty}) \|v\|_{H^m}  \|Dv\|_{H^m} \n \\
 && \leq \frac12 \|Dv\|_{H^m} ^2 + C(\|v\|_{L^\infty}^2 +\|\nabla v\|_{L^\infty}^2) \|v\|_{H^m} ^2,
\eqn
and therefore
\bq
\lefteqn{\| v(t)\|_{H^m} ^2 +\int_0 ^t \|Dv(s)\|_{H^m} ^2 ds \leq\| v_0\|_{H^m} ^2\exp\left(C\int_0 ^t
(\|v\|_{L^\infty}^2 +\|\nabla v\|_{L^\infty}^2)ds\right)}\n \\
&&\leq \| v_0\|_{H^m} ^2\exp\left\{C\int_0 ^t
(\|v\|_{L^2}^2 +\|D^3 v\|_{L^2}^2)ds\right\}\n \\
&&\leq \| v_0\|_{H^m} ^2\exp \left\{t \|v_0\|_{L^2}^2 \exp (C\|\Lambda ^{\frac12 }
 v_0\|_{L^2}^2)\right\}\times\n \\
  &&\times\exp\left[ \|\Lambda v_0\|_{L^2} ^2\exp \left\{Ct\|v_0\|_{L^2}^2 \exp (C\|\Lambda ^{\frac12 }
 v_0\|_{L^2}^2)  \|\Lambda v_0\|_{L^2} ^2\exp (C\|\Lambda ^{\frac12 }
 v_0\|_{L^2}^2)\right\}\right],\n \\
\eq
where we used the estimates (\ref{est1}) and (\ref{est3}). $\square$
 $$ \mbox{\bf Acknowledgements} $$
This work was supported partially by  NRF Grant no. 2006-0093854 and  2009-0083521.

\end{document}